\documentclass{article}

\usepackage{amsmath,amssymb,dsfont,amscd,theorem,graphicx}
\usepackage[all]{xy}
\usepackage{color}

\long\def\footinfo#1{
\begingroup\def\thefootnote{\fnsymbol{footnote}}\footnote[0]{\hskip-18pt #1}\endgroup}

\def\suported#1{#1.}
\def\suport#1{Supported in part by UI\&D {\sl Matem\'{a}tica e aplica\c{c}\~{o}es} of University of
Aveiro, through Program POCTI of FCT co-financed by the European Community fund FEDER.}
\def\ZOOOsubjclass#1{2000 {\em Mathematics Subject Classification}. #1}
\def\keywords#1{ {\em Key words and phrases}: #1.}

\def\oop#1{}
\def\hide#1{}

\def\pmeta#1{M{(#1)}}
\def\pmetaem#1{{\em M}{({\em #1})}}

\newtheorem{thm}{Theorem}
\newtheorem{cor}[thm]{Corollary}
\newtheorem{lem}[thm]{Lemma}
\theorembodyfont{\upshape}
\newtheorem{rem}{Remark}

\newenvironment{proof}{\noindent\textsc{Proof.}}{$\blacksquare$}

\newcounter{diag}
\newcounter{tab}

\begin{document}

\title{Strong double coverings of groups}
\author{Ana Breda, Antonio Breda d'Azevedo, Domenico Catalano}
\maketitle

\footinfo{\ZOOOsubjclass{20E22}} \footinfo{\keywords{Group extensions, Metacyclic groups, simple groups, platonic
groups, strong coverings}} \footinfo{\suported{Supported in part by UI\&D {\sl Matem\'{a}tica e aplica\c{c}\~{o}es} of
University of Aveiro, through Program POCTI of FCT co-financed by the European Community fund FEDER}}

\begin{abstract}
By a covering of a group $G$ we mean an epimorphism from a group $\hat G$ to $G$. Introducing the notion of strong
covering as a covering \mbox{$\pi:\hat G\rightarrow G$} such that every automorphism of $G$ is a projection via $\pi$
of an automorphism of $\hat G$, the main aim of this paper is to characterise double coverings which are strong. This
is done in details for metacyclic groups, rotary platonic groups and some finite simple groups.
\end{abstract}

\section{Preliminaries}

\hide{In this paper we develop some group theoretical tools useful in hypermap theory to count and
describe the rotary double coverings of a given rotary hypermap. (For an introduction to hypermaps
see \cite{B1,CS,IS,LJ,GJ}). In fact, the rotary double coverings of a given hypermap are
completely described by the double coverings of its automorphism group. Using this, and the
results obtained in this paper, we intend to give examples of application to hypermaps in a next
paper, limiting us here to the group theoretical aspects including examples (\S\ref{DCMG},
\S\ref{DCPG} and \S\ref{DCSG}).}

In map/hypermap theory, or in polytopes theory, it is sometimes desirable to know not only the
coverings $\hat G$ of a given group $G$ but also whether the automorphism group $Aut(\hat G)$ of
the covering also covers the automorphism group of $G$. Such coverings will be called ``strong
coverings''. In this paper we address this problem to double coverings. The characterisation of
strong double coverings will be done in details for metacyclic groups, rotary platonic groups and
some finite simple groups.

\medskip

For simplicity when dealing with presentations of groups we need to establish some conventions. We will generally see
the set of generators and the set of relators of a given group presentation $\langle x_1,...,x_n\mid
r_1,...,r_m\rangle$ as ordered sets and, as such, as tuples. Hence, setting $X=(x_1,...,x_n)$ and $R=(r_1,...,r_m)$ we
will write $\langle X\mid R\rangle$ to mean $\langle x_1,...,x_n\mid r_1,...,r_m\rangle$ and call $n$ the \emph{rank}
of the presentation $\langle X\mid R\rangle$. Moreover by $r_k(X)$ and $R(X)$ we mean $r_k(x_1,...,x_n)$ and
$(r_1(X),...,r_m(X))$, respectively. Seeing the set of generators and the set of relators of a group presentation as
ordered sets in the form of tuples, we can bring the cartesian product notations to presentations. For instance, if
$R=(r_1,...,r_m)$ and $S=(s_1,...,s_m)$ are ordered set of words on $x_1,...,x_n$ then $R=S$ has the usual ``cartesian"
meaning $r_k=s_k$, $k=1,...,m$. For convenience, especially when dealing with equations, we may use relations instead
of relators in presentations. We reserve the notation $[a,b]$ for the commutator $a^{-1}b^{-1}ab$ of elements $a,b$ of
a group $G$ and the natural simplification $[a,X]$ to mean $[a,x_1],...,[a,x_n]$ for the $n$-tuple $X=(x_1,...,x_n)$ of
elements of $G$. Finally, given two presentations $P_1$ and $P_2$ we write $P_1\sim P_2$ if they are presentations of
the same group.

\medskip

Let $\Gamma$ be a group with presentation $\langle X\mid W\rangle$ where $X=(x_1,...,x_n)$ and $W$ is an ordered set of
words on $x_1,...,x_n$. By a \emph{$\Gamma$-base} of a finite group $G$ we mean a $n$-tuple $(a_1,...,a_n)$ of elements
of $G$ for which there is an epimorphism from $\Gamma$ to $G$ mapping $x_i$ to $a_i$\,. If $\Gamma$ is the free group
of rank $n$ we say \emph{$n$-base} instead of $\Gamma$-base. Hence a $n$-base of a group $G$ is just a $n$-tuple
$(a_1,...,a_n)$ of elements of $G$ such that $\{a_1,....,a_n\}$ generates $G$. We will denote by $B_n(G)$ the set of
$n$-bases of $G$ and say that $G$ is \emph{$n$-generated} if $B_n(G)\neq\varnothing$\,.

\medskip

Given a presentation $P=\langle X\mid R\rangle$ of rank $n$ of a group $G$ we denote by $S_P$ the set of $n$-bases
$(a_1,...,a_n)$ of $G$ satisfying $R(a_1,...,a_n)=1$\,. The elements of $S_P$ will be called \emph{presentation
$n$-tuples of $P$}. Obviously every $n$-base of $G$ belongs to a presentation set $S_P$ for some presentation $P$ of
$G$ of rank $n$.

\begin{thm}\label{SamePSet}
$(a_1,...,a_n),(b_1,...,b_n)\in B_n(G)$ belong to the same presentation set if and only if there is an automorphism of
$\,G$ mapping $a_i$ to $b_i$.
\end{thm}

\begin{proof}
If $(a_1,...,a_n),(b_1,...,b_n)\in B_n(G)$ belong to the same presentation set $S_P$ then by the Substitution Test the
functions $x_i\mapsto a_i$ and $x_i\mapsto b_i$ extend to isomorphisms $\alpha$ and $\beta$ from $P$ to $G$,
respectively. Then $\alpha^{-1}\beta$ is an automorphism of $G$ mapping $a_i$ to $b_i$. The converse is
straightforward.
\end{proof}

\begin{cor}\label{Aut=SP}
$|Aut(G)|=|S_P|$ for every presentation $P$ of finite rank of the group $G$.
\end{cor}

Each function $f:A\rightarrow B$ extends, in a natural way, to a function
\mbox{$f^*:A^n\rightarrow B^n$} defined by $(a_1,...,a_n)\mapsto (a_1f,...,a_nf)$. As usual, we
will denote $f^*$ also by $f$, provided no confusion arises from such simplification. Hence, for
each $n$-tuple $b\in B^n$, $bf^{-1}$ will denote the set $\{a\in A^n\mid af^*=b\}$.

\section{Presentations of double coverings}

Let $p\in\mathds{N}$. An epimorphism $\pi:\hat G\rightarrow G$ will be called a \emph{$p$-covering (of $G$)} if it has
kernel of size $p$ and a \emph{central covering (of $G$)} if its kernel is in the center of $\hat G$.

\begin{rem}\label{CovPair}
A $p$-covering of $G$ is determined by a pair $(\hat G,N)$ consisting of a group $\hat G$ and a normal subgroup $N$ of
$\hat G$ of cardinality $p\,$ such that $\hat G/N$ is isomorphic to $G$; best known in the literature as an
\emph{extension} of $G$ by $N$.
\end{rem}

\noindent The proof of the following Lemma is straightforward.

\begin{lem}\label{Double=Central}
Every double covering ($2$-covering) is a central covering.
\end{lem}

From now on, let $P=\langle X\mid R\rangle$ be a presentation of a group $G$, where $X=(x_1,...,x_n)$ and
$R=(r_1,...,r_m)$ and $\pi:\hat G\rightarrow G$ be a double covering. Since $\langle i\mid i^2\rangle$ is a
presentation of $Kern(\pi)\cong C_2$, which lies in the center of $\hat G$ (Lemma \ref{Double=Central}), applying
Johnson \cite{DJ}, Chapter 10, to this case we get the following result.

\begin{lem}\label{Johnson}
For every presentation $n$-tuple $a\in S_P$ and every $\hat a=(\hat a_1,...,\hat a_n)\in a\pi^{-1}$ there is
$J\in\{1,i\}^m=C_2^m$ such that
\begin{equation}\label{pres}
P_J=\langle X,i\mid R=J,i^2=1,[i,X]=1\rangle
\end{equation}
is a presentation of $\hat G$ with $(\hat a,i)=(\hat a_1,...,\hat a_n,i)\in S_{P_J}$
\end{lem}

\noindent For simplicity we will omit in the expression \eqref{pres} the obvious relations $i^2=1$ and $[i,X]=1$ that
the central involution $i$ naturally satisfies, and write
\[
P_J=\langle X,i\mid R=J\,\rangle.
\]
If $J=1$ then $P_J$ is a presentation of $G\times C_2$ of rank $n+1$. If $J\neq 1$ then
there is $k\in\{1,...,m\}$ such that $r_k=i$. Thus $i\in\langle x_1,\dots,x_n\rangle$ and
therefore we write
\[
P_J=\langle X\mid R=J\,\rangle
\]
and regard $P_J$ as a presentation of rank $n$ with $\hat a\in S_{P_J}$.

\begin{rem}
For every $J\in C_2^m$ we have $P_J/\langle i\rangle\cong G$ and therefore $(P_J,\langle i\rangle)$ is, according to
Remark \ref{CovPair}, an extension of $G$ by the subgroup of $P_J$ generated by $i$. Two cases may occur: either $i$ is
in the normal closure $K$ generated by the relators of $P_J$ in the free group $F(X,i)$ or not. In the first case $P_J$
is a presentation of $G$. In the second case $i\not\in K$ and $(P_J,\langle i\rangle)$ is an extension of $G$ by
$C_2\cong\langle i\rangle$ i.e. $P_J$ is a double covering of $G$. \oop{In this case, $P_J$ or the group $\hat G$ with
presentation $P_J$ will be called a double covering of $G$ since the double covering $\pi:\hat G\rightarrow G$ is then
implicitly given by $P_J$ (and $\langle i\rangle$).} We call $P_{(i,\dots,i)}$ the \emph{binary} of $G$, and denote it
by $\widetilde{G}$, if it is a double covering of $G$. Not every such presentation gives always a binary group, there
are cases in which $P_{(i,\dots,i)}$ collapses. For example, $G=PSL(2,2^f)$, for $f>2$, and some metacyclic groups $G$
(see Table 2 where $[P_{(i,i,i)}]=[P_{(1,1,i)}]$ in \#\textbf{1}, \textbf{5} and $[P_{(i,i,i)}]=[P_{(1,i,i)}]$ in
\#\textbf{7}) the binary presentation collapses to a presentation of $G$. In the case of $G=D_n,A_4,S_4$ and $A_5$
the binary presentation $P_{(i,i,i)}$ coincides with the usual binary group $\widetilde{G}$ (see \cite{La}).
\end{rem}

Let $I=(i_1,...,i_n)$ be an element of the elementary 2-group $Kern(\pi)^n\cong C_2^n$. Seeing $X$ as an element of
$\hat G^n$ we may write $IX=(i_1x_1,...,i_nx_n)$. Since $i$ is in the centre of $\hat G$ we have
\[
r_k(IX)=r_k(i_1x_1,...,i_nx_n)=r_k(i_1,...,i_n)\,r_k(x_1,...,x_n)=r_k(I)\,r_k(X)
\]
and therefore $R(IX)=R(I)\,R(X)$ in the direct product $F(X,i)^m$.

\begin{lem}\label{congpres}
Let $J\in C_2^m$ and $P_J$ be a presentation of $\hat G$ (a double covering of $G$). For any $I\in C_2^n$, $P_{R(I)J}$
is another presentation of $\hat G$; that is, $P_{R(I)J}\sim P_J$.
\end{lem}

\begin{proof}
By changing generators $(X,i)$ to $(Y,i)=(IX,i)$ we get $P_J=\langle X,i\mid R=J\,\rangle\sim\langle Y,i\mid
R=R(I)J\rangle=P_{R(I)J}$. \oop{(Note that this does not produce extra relations since $Y=I^2Y=Y$ and
$R(Y)=R(I)\,R(X)=R(I)J$).}
\end{proof}

\bigskip

\noindent Let $R(C_2^n)=\{R(I)\mid I\in C_2^n\}\subset C_2^m$. It is easy to check that $R(C_2^n)$ is a normal subgroup
of $C_2^m$. The following corollaries are easy consequence of Lemma \ref{congpres}.

\begin{cor}\label{NPC}
Given a group $\,G$ with presentation $\langle X\mid R\,\rangle$\oop{, where $X=(x_1,...,x_n)$ and $R=(r_1,...,r_m)$,}
the number of presentation classes is given by the index of $R(C_2^n)$ in $C_2^m$.
\end{cor}

\begin{cor}\label{NDC}
Given a group $\,G$ with presentation $\langle X\mid R\,\rangle$\oop{, where $X=(x_1,...,x_n)$ and $R=(r_1,...,r_m)$,}
the number of (non-isomorphic) double coverings of $\,G$ is less or equal than the index $|C_2^m:R(C_2^n)|$.
\end{cor}

\noindent Equality is not always reached as the examples in \S\ref{DCMG} will show. However, according to Lemma
\ref{congpres}, we may cluster the $2^m$ presentations $P_J$, $J\in C_2^m$, in equivalence classes given by the
equivalence relation
\[P_J\approx P_L\,\, \Leftrightarrow \,\,L\in R(C_2^n)J\,.\]
Each equivalence class $[P_J]$ contains $|R(C_2^n)|$ presentations of the same group $\hat G$ and therefore $[P_J]$
will be called a \emph{presentation class of $\,\hat G$}. As remarked above, if $P_J$ is not a presentation of a double
covering of $G$ then it is a presentation of $G$. We observe that this can not happen if $P_J\in [P_1]$ since $P_1$
(and therefore $P_J$) is a presentation of the direct product $\hat G=G\times C_2$. Two distinct presentation classes
may be presentation classes of the same double covering $\hat G$. However, in the particular case of $\hat G=G\times
C_2$ this can not happen as the following theorem shows.

\begin{thm}\label{UniquePClass}
Let $P$ be a presentation of the group $G$. Then $\,G\times C_2$ has only one presentation class, namely $[P_1]$.
\end{thm}

\begin{proof}
Let $\pi:G\times C_2\rightarrow G$ be the canonical projection. Suppose that $J\in C_2^m$ gives rise to a presentation
$P_J\not\in [P_1]$ of $G\times C_2$. In this case each presentation in $[P_J]$ has rank $n$. By Lemma \ref{Johnson}
there are $a\in S_P$ and $\hat a\in a\pi^{-1}$ such that $\hat a\in S_{P_J}$. Then by Lemma \ref{Johnson} and Lemma
\ref{congpres}, for every $\hat x\in a\pi^{-1}$ there is $P_L\in [P_J]$ such that $\hat x\in S_{P_L}$. This shows that
$a\pi^{-1}\subset B_n(G\times C_2)$, which is a contradiction since $\hat x=((a_1,1),...,(a_n,1))\in a\pi^{-1}\setminus
B_n(G\times C_2)$. Hence $P_J\in [P_1]$.
\end{proof}

\bigskip

For practical reasons, it is convenient to simplify, whenever possible, the presentation $P$ of $G$ in order to
emphasize the number of cosets of $R(C_2^n)$ in $C_2^m$. An obvious simplification is given by choosing an
``efficient'' presentation of $G$, i.e. a presentation with maximal deficiency, where the deficiency of a presentation
is the non-positive integer given by the difference between the rank and the number of relators of the presentation.
Another possibility is given by presentations of ``simple type'' defined below.

\medskip

Let $w$ be an element of the free group $F(n)=F(x_1,...,x_n)$. By the \emph{$i^{th}$-exponent-sum} of $w$ we mean the
image $w_i$ of $w$ under the epimorphism \mbox{$F(n)\rightarrow\mathds{Z},\ x_j\mapsto\delta_{ij}$,} $j=1,...,n$, where
$i\in\{1,...,n\}$ and $\delta_{ij}=1$ if $i=j$ and $0$ otherwise.\oop{ The index $i$ of $w_i$ will be called the
\emph{exponent-sum index}\,.}

\medskip

If a relator $w=r_k$ of $P$ has even exponent-sum $w_i$, for all $i$, we say that $w$ is \emph{even}, otherwise we say
that $w$ is \emph{odd}. Let $E_P$ and $O_P$ be respectively the subset of even and odd relators of $P$. Of course
$E_P\cup O_P$ is the set of all relators of $P$. We say that $P$ is of \emph{simple type} if:
\begin{itemize}
\item Each $w\in O_P$ gives rise to only one odd $i^{th}$-exponent-sum $w_i$; \oop{that is, for each odd relator $w$ of $P$
there is only one exponent-sum index $i$ such that $w_i$ is odd,}
\item For each $i\in \{1,\dots,n\}$ there is at most one $w\in O_P$ such that $w_i$ is odd.
\end{itemize}
In this case we call the number $d=|O_P|$ the \emph{degree} of the presentation $P$.

\begin{thm}\label{SimpleType}
If $P$ is a presentation of simple type of $G$ of rank $n$ with $m$ relators and degree $d$ then $d\leqslant n$ and the
index of $R(C_2^n)$ in $C_2^m$ is $2^{|E_P|}=2^{m-d}$.
\end{thm}

\begin{proof}
That $d=|O_P|\leqslant |X|=n$ follows from the definition. Suppose, without loss of generality, that
$O_P=(r_1,...,r_d)$ and $E_P=(r_{d+1},...,r_m)$. Then $R(C_2^n)=O_P(C_2^n)\times E_P(C_2^n)$ is a normal subgroup of
$C_2^m=C_2^{|O_P|+|E_P|}=C_2^{|O_P|}\times C_2^{|E_P|}$. Now, since $O_P(C_2^n)\cong C_2^{|O_P|}$ and
$E_P(C_2^n)=\{(1,...,1)\}< C_2^{|E_P|}$, we have that $|C_2^m:R(C_2^n)|=|C_2^{|E_P|}:E_P(C_2^n)|=2^{|E_P|}=2^{m-d}$.
\end{proof}

\medskip

Let $P$ be a presentation of simple type of $G$ of degree $d$. By Theorem \ref{SimpleType} and Corollary \ref{NPC}, $P$
gives rise to $2^{|E_P|}$ presentation classes. Analysing the structure of the elements of each presentation class we
get the following useful result for presentation of simple type.

\begin{cor}\label{CorSimpleType}
Each presentation class $[P_J]$ contains $2^d$ presentations and is represented by a unique presentation of the form
$P_L=\langle X,i\mid O_P=1,E_P=L\rangle$
with $L=(i_{d+1},...,i_m)\in\{1,i\}^{m-d}$, where $m$ is the number of relators of $P$.
\end{cor}

\section{Strong double coverings}

A $p$-covering $\pi:\hat G\rightarrow G$ induces an equivalence relation $\sim_\pi$ on $\hat G$ defined by $\hat
g_1\sim_\pi\hat g_2\Leftrightarrow\hat g_1\pi=\hat g_2\pi$. The equivalence classes are the cosets $[\hat g]_\pi=\hat
g\,Kern(\pi)$ and the class partition is $$\hat G/\pi=\hat G/\sim_\pi=\hat G/Kern(\pi)\,\cong G\,.$$ If $\psi$ is an
automorphism of $\hat G$ satisfying $Kern(\pi)\,\psi\subset Kern(\pi)$ then
\[
\psi^\pi:\hat G/\pi\rightarrow\hat G/\pi,\ [\hat g]_\pi\mapsto [\hat g\psi]_\pi
\]
is a well-defined function satisfying $\hat g\pi\,\psi^\pi=\hat g\psi\pi$ for every $\hat g\in \hat G$ i.e. the
following diagram commutes
\[
\begin{CD}
\hat G @>{\psi}>> \hat G\\[0pt]
@V{\pi}VV@VV{\pi}V\\[0pt]
G@>{\psi^\pi}>> G
\end{CD}
\]

\medskip

\noindent It is straightforward to show that $\psi^\pi$ is an automorphism. Since
$\psi_1^\pi\,\psi_2^\pi=(\psi_1\psi_2)^\pi$ for every $\psi_1,\psi_2\in Aut(\hat G)$ we have the following statement.

\begin{lem}\label{CharKern}
If $Ker(\pi)$ is a characteristic subgroup of $\hat G$, then \mbox{$\zeta^\pi:Aut(\hat G)\rightarrow Aut(G),\
\psi\mapsto\psi^\pi$}, is a well-defined homomorphism.
\end{lem}

When well-defined, $\zeta^\pi$ maps inner automorphisms of $\hat G$ to inner automorphism of $G$, but in general it may
be not onto. If \mbox{$\pi:\hat G\rightarrow G$} is a $p$-covering with characteristic kernel and
\mbox{$\zeta^\pi:Aut(\hat G)\rightarrow Aut(G)$} is an epimorphism we say that $\pi$ is a \emph{strong covering}, or a
\emph{$q$-strong $p$-covering} if $\zeta^\pi$ is a $q$-covering.

\medskip

Let $\pi:\hat G\rightarrow G$ be a $p$-covering and let $\hat a\in S_{\hat P}$ where $\hat P$ is a presentation of
$\hat G$. Further let $P$ be a presentation of $G$ with $a=\hat a\pi\in S_P$.

\begin{thm}\label{sizeker}
If $Kern(\pi)$ is a characteristic subgroup of $\hat G$, then
\begin{itemize}
\item[\emph{(a)}] $\ \hat b\pi\in S_P$ for every $\hat b\in S_{\hat P}$\,;
\item[\emph{(b)}] $\ |Kern(\zeta^\pi)|=|a\pi^{-1}\cap S_{\hat P}|$\,;
\item[\emph{(c)}] $\ \zeta^\pi$ is an epimorphism if and only if $\,b\pi^{-1}\cap S_{\hat P}\neq\varnothing\,$ for every $b\in S_P$\,.
\end{itemize}
\end{thm}

\begin{proof}
(a) Let $\hat b\in S_{\hat P}$ and, according to Theorem \ref{SamePSet}, let $\psi\in Aut(\hat G)$ such that $\hat
a\psi=\hat b$. Then\oop{, from $a\in S_P$,} $\hat b\pi=\hat a\psi\pi=\hat a\pi\psi^\pi=a\psi^\pi$ and by Theorem
\ref{SamePSet}, $\hat b\pi\in S_P$ since $\psi^\pi=\psi\zeta^\pi\in Aut(G)$.\\[3pt]
(b) $\psi\in Kern(\zeta^\pi)\Leftrightarrow a\psi^\pi=a\Leftrightarrow\hat a\pi\psi^\pi=a\Leftrightarrow\hat
a\psi\pi=a\Leftrightarrow\hat a\psi\in a\pi^{-1}$. As $\hat a\in S_{\hat P}$ the statement follows from Theorem
\ref{SamePSet} and the fact that for every $\psi_1,\psi_2\in Aut(\hat G)$, $\hat a\psi_1=\hat
a\psi_2\Leftrightarrow\psi_1=\psi_2$ (Corollary \ref{Aut=SP}).\\[3pt]
(c) $\zeta^\pi$ is onto if and only if for every $\phi\in Aut(G)$ there is $\psi\in Aut(\hat G)$ such that
$\psi^\pi=\phi$. According to Theorem \ref{SamePSet} this is equivalent to
\[
\forall\ b\in S_P\,,\ \exists\ \hat b\in S_{\hat P}\,,\ \hat b\pi=b
\]
proving that $\zeta^\pi$ is onto if and only if $\,b\pi^{-1}\cap S_{\hat P}\neq\varnothing\,$ for every $b\in S_P$\,.
\end{proof}

\begin{thm}\label{strong}
Let $\pi:\hat G\rightarrow G$ be a double covering with characteristic kernel and let $P$ be a presentation of $\,G$.
If $\,\hat G$ has a presentation $P_J$ with $J\neq 1$ then $\pi$ is a strong double covering of $\,G$ if and only if
$\,[P_J]$ is the only presentation class of $\,\hat G$.
\end{thm}

\begin{proof}
Let $[P_L]$ be a presentation class of $\hat G$ with $L\in C_2^m$, where $m$ is the number of relators of the
presentation $P$. Then $[P_L]\neq \{P_1\}$ since otherwise $\hat G\cong G\times C_2$ and Theorem \ref{UniquePClass}
will contradict the hypothesis that there is a presentation $P_J$ of $\hat G$ with $J\neq 1$. Hence we can assume
$L\neq 1$ and regard
$P_L$ as a presentation of rank $n$, where $n$ is the rank of the presentation $P$.\\[3pt]
Suppose that $\pi$ is a strong double covering of $G$ and let $\hat b\in S_{P_L}$. Then $b=\hat b\pi\in S_P$ and
$b\pi^{-1}=\{I\hat b\mid I\in C_2^n\}$. As $\zeta^\pi$ is an epimorphism then, by Theorem \ref{sizeker} (c), there is
$I\in C_2^n$ such that $I\hat b\in S_{P_J}$. This implies $J=R(I\hat b)=R(I)R(\hat b)=R(I)L$. Hence $[P_J]=[P_L]$.\\[3pt]
Reciprocally, suppose that $[P_J]$ is the only presentation class of $\hat G$. Let $b\in S_P$ and $\hat b\in
b\pi^{-1}$. Then $R(\hat b)=L$ for some $L\in C_2^m$. As $[P_J]$ is the only presentation class of $\hat G$ then
$L=R(I)J$ for some $I\in C_2^n$. Hence $J=R(I)L=R(I)R(\hat b)=R(I\hat b)$ and therefore $I\hat b\in b\pi^{-1}\cap
S_{P_J}$. By Theorem \ref{sizeker} (c) we conclude that $\zeta^\pi$ is an epimorphism and therefore $\pi$ is a strong
double covering of $G$.
\end{proof}

\begin{cor}\label{q=SizeKrho}
Let $\pi:\hat G\rightarrow G$ be a $q$-strong double covering of $\,G$ and let $P=\langle X\mid R\rangle$ be a
presentation of rank $\,n$ of $\,G$ with $\,m$ relators. If $\,\hat G$ has a presentation $\,P_J$ for some $J\neq 1$
then $\,q=|Ker(\rho)|$\,, where \mbox{$\,\rho:C_2^n\rightarrow C_2^m,\ I\mapsto R(I)$\,.}
\end{cor}

\begin{proof}
Let $\pi:\hat G\rightarrow G$ be a strong double covering of $G$ and let $\hat a\in S_{\hat P}$ where $\hat P=P_J$ is a
presentation of $\hat G$ with $J\neq 1$. Then $\hat P$ is a presentation of rank $n$ of $\hat G$ and, by Theorem
\ref{sizeker} (a), $a=\hat a\pi\in S_P$. According to Lemma \ref{Johnson}, for every $\hat b\in a\pi^{-1}$ there is
$L\in C_2^m$ such that $P_L$ is a presentation of $\hat G$ with $\hat b\in S_{P_L}$. By Theorem \ref{strong}, $P_L\in
[P_J]$, i.e. $L=R(I)J$ for some $I\in C_2^m$. Hence
\[
a\pi^{-1}\cap S_{\hat P}=\{I\hat a\mid I\in C_2^n\}\cap S_{P_J}=\{I\hat a\mid I\in C_2^n,\ R(I)=1\}
\]
and the statement follows from Theorem \ref{sizeker} (b).
\end{proof}

\medskip

\noindent Under the assumptions of Corollary \ref{q=SizeKrho}, if $P$ is a presentation of simple type we may
deduce from Theorem \ref{SimpleType} the following corollary.

\begin{cor}\label{q-simpletype}
Let $\,\pi:\hat G\rightarrow G$ be a $q$-strong double covering of $\,G$. If $P$ is a
presentation of simple type of degree $d$ (and rank $n$) of $\,G$ and $\,\hat G$ has presentation $\,P_J$ for some $J\neq 1$,
then $q=2^{n-d}$\,.
\end{cor}
\oop{\begin{proof} In fact, $\rho:C_2^n\rightarrow Im(\rho)=R(C_2^n)$ and so $|Im(\rho)|=|Im(\rho):1|=|C_2^n:Ker(\rho)|$. In other
words, $\frac{2^n}{|Ker(\rho)|}=|Im(\rho)|\Leftrightarrow |Ker(\rho)|=\frac{2^n}{|Im(\rho)|}$. By Theorem
\ref{SimpleType}, $|C_2^m:Im(\rho)|=2^{m-d}$, i.e. $|Im(\rho)|=2^d$. Hence $|Ker(\rho)|=\frac{2^n}{2^d}=2^{n-d}$.\end{proof}}

\subsection*{The direct product \textbf{\emph{G}}$\,\times$\textbf{\emph{C}}$_\mathbf{2}$}

Remark that if $\pi:\hat G\rightarrow G$ is a double covering and $G$ has no central involution, then $Kern(\pi)$ is a
characteristic subgroup of $\hat G$. The following Theorem shows that the converse is true for $\hat G=G\times C_2$.

\begin{thm}\label{KerCharGxC2}
The kernel of the canonical epimorphism $\pi:G\times C_2\rightarrow G$ is a characteristic subgroup of $\,\hat G$ if
and only if $\,G$ has no central involution.
\end{thm}

\begin{proof}
Let $P=\langle X\mid R\rangle$ be a presentation of $G$. If $G$ has a central involution $w$ then, since $P_1=\langle
X,i\mid R,i^2,[i,X]\rangle$ is a presentation of $G\times C_2$ and $(X,wi)\in S_{P_1}$, by Theorem \ref{SamePSet},
there is an automorphism (outer) fixing $X$ and sending $i$ to $wi$ i.e. $Ker(\pi)=\{1,i\}$ is not a characteristic
subgroup of $G\times C_2$.
\end{proof}

\medskip

For every group $G$ and for every $\varphi\in Aut(G)$ we note that $\varphi\times id\in Aut(G\times C_2)$. Hence
$\zeta^\pi$ is onto, if well-defined. Combining this with Lemma \ref{CharKern} we get the following corollary of
Theorem \ref{KerCharGxC2}.

\begin{cor}\label{StrongGxC2}
The canonical epimorphism $\pi:G\times C_2\rightarrow G$ is a strong double covering of $G$ if and only if $G$ has no
central involutions.
\end{cor}

We conclude this section with some results about the special case when $G$ is a simple group without central
involutions, i.e. not isomorphic to $C_2$.

\begin{thm}\label{PropNormSubgroups}
If $G$ is a simple group not isomorphic to $C_2$ then $G\times\{1\}$ and $\{1\}\times C_2$ are the only proper normal
subgroups of $G\times C_2$.
\end{thm}

\begin{proof}
Let $\pi$ be the canonical epimorphism from $G\times C_2$ to the simple group $G\not\cong C_2$ and let $N$ be a proper
normal subgroup of $G\times C_2$. Then $N\pi$ is a normal subgroup of $G$ and the simplicity of $G$ implies
$N\pi=\{1\}$ or $G$. If $N\pi=\{1\}$ then $N=\{1\}\times C_2$. If $N\pi=G$ then $N$ has index 2 in $G\times C_2$. As
$N\cap (G\times\{1\})$ is a normal subgroup of the simple group $G\times\{1\}$ we conclude that $N=G\times\{1\}$.
\end{proof}

\begin{cor}
If $G$ is a simple group not isomorphic to $C_2$ then \mbox{$|Aut(G\times C_2)|=|Aut(G)|$.}
\end{cor}
\oop{
\begin{proof}
This occur from $|Aut(G\times C_2)|=|S_{P_1}|=|S_P|=|Aut(G)|$. In fact, $(X,i), (X',i')\in S_{P_1}$ implies one of $i$
or $i'$ belongs to $G=G\times\{1\}$ and the other to $C_2=\{1\}\times C_2$. Then $\langle i,i'\rangle< C_2\times C_2$
and $\langle i,i'\rangle\lhd G\times C_2$. Since $G$ is not $C_2$ then by Theorem \ref{PropNormSubgroups} $\langle
i,i'\rangle=1\times C_2$ and so $i=i'$. Hence counting presentation tuples $(X,i)$ in $S_{P_1}$ is the same as counting
presentation tuples $X$ in $S_{P}$.
\end{proof} \medskip
}

\noindent Since a simple group not isomorphic to $C_2$ has no central involution we have the following corollary.

\begin{cor}\label{SimpleStrong}
If $\,G$ is a simple group not isomorphic to $C_2$, then the canonical epimorphism $G\times C_2\rightarrow G$ is a
1-strong double covering.
\end{cor}

\section{Double coverings of metacyclic groups}\label{DCMG}

A group $G$ is called \emph{metacyclic} if it has a normal subgroup $N$ such that both $N$ and $G/N$ are cyclic. It is
widely known \cite{DJ} that the group $\,\pmeta{m,n,r,s}=\langle x,y\mid x^m=1, y^n=x^r, x^y=x^s\rangle\,$ with parameters $m,n,r,s\in\mathds{N}$
is a metacyclic group of order less or equal to $mn$ with equality if and only if the parameters
satisfy
\begin{equation}\label{metacond}
s^n\equiv 1\mod m\qquad\text{and}\qquad rs\equiv r\mod m\,.
\end{equation}
Moreover, every metacyclic group has a presentation of the form
$\pmeta{m,n,r,s}$ with $m,n,r,s$ satisfying \eqref{metacond}, where we can obviously assume
$r,s\leqslant m$. Cyclic groups, dihedral groups and dicyclic groups are just particular cases of metacyclic groups.
So, as an introduction to the general case, we will first handle the cyclic, dihedral and dicyclic groups, in this
order. But before that, we give an observation, compute the central involutions and give the presentations pairs for
$\pmeta{m,n,r,s}$.

\subsection{Observation}

Unless otherwise clearly specified, $M(m,n,r,s)$ is a metacyclic group of order $mn$.
Directly from the presentation $\pmeta{m,n,r,s}$ we get that for all $\alpha,\beta\in\mathds{N}$
\begin{equation}\label{perm(x,y)}
x^\alpha y^\beta =y^\beta x^{\alpha s^\beta}\,.
\end{equation}
Thus $\ \pmeta{m,n,r,s}=\{y^px^q\mid 0\leqslant p<n,0\leqslant q<m\}\ $ and so, as a set,
$\pmeta{m,n,r,s}=\mathds{Z}_n\times\mathds{Z}_m$. From \eqref{perm(x,y)} we get the following power formula
\begin{equation}\label{powers}
\left( y^px^q\right)^k=y^{p\,k}\,x^{q\,\sigma(s^p,k)}\,,
\end{equation}
where $\ \sigma:\mathds{N}\times\mathds{N}\rightarrow\mathds{N},\
(h,k)\mapsto\sum\limits_{i=0}^{k-1}h^i=\left\{\begin{array}{ll} \ k & \text{if }\ h=1\\ \frac{\ h^k-1}{h-1} &
\text{otherwise}\end{array}\right.$\,.

\subsection{The central involutions of \pmetaem{m,n,r,s}}\label{CentralInv}

By \eqref{powers} if $m$ and $n$ are both odd, then $G=\pmeta{m,n,r,s}$ has no involution at all. If $m$ is even, then
$x^{m/2}$ is an involution. According to \eqref{metacond}, $m$ even implies $s$ odd and from \eqref{perm(x,y)} we get
\[
x^{\frac{m}{2}}y=yx^{\frac{m}{2}s}=yx^{\frac{m}{2}}\,,
\]
which shows that $x^{m/2}$ belongs to the center. If in addition $n$ is odd then $x^{m/2}$ is the unique involution in
$G$. If $n$ is even and $p\neq 0$ then $y^px^q$ is an involution if and only if $p=n/2$ and
\[
q(s^{\frac{n}{2}}+1)+r\equiv 0\mod m\,.
\]
Moreover, $y^{n/2}x^q$ is in the center if and only if
\[
\left\{\begin{array}{l}
(y^{n/2}x^q)\,x=x\,(y^{n/2}x^q)\\
(y^{n/2}x^q)\,y=y\,(y^{n/2}x^q)
\end{array}\right.\ \Leftrightarrow\
\left\{\begin{array}{l}
x=x^{s^{\frac{n}{2}}}\\
x^{qs}=x^q
\end{array}\right.\ \Leftrightarrow\
\left\{\begin{array}{ll} s^{\frac{n}{2}}\equiv 1 & \hspace*{-2ex}\mod m\\
q(s-1)\equiv 0 & \hspace*{-2ex}\mod m \end{array}\right.
\]
Hence $y^{n/2}x^q$ is a central involution if and only if
\[
\left\{\begin{array}{ll} q(s^{\frac{n}{2}}+1)+r\equiv 0 & \hspace*{-2ex}\mod m\\
s^{\frac{n}{2}}\equiv 1&\hspace*{-2ex}\mod m\\
q(s-1)\equiv 0 &\hspace*{-2ex}\mod m
\end{array}\right.\quad\Leftrightarrow\quad
\left\{\begin{array}{lll} 2q+r\equiv 0 & \hspace*{-2ex}\mod m & (\star)\\
s^{\frac{n}{2}}\equiv 1&\hspace*{-2ex}\mod m & (\ast)\\
q(s-1)\equiv 0 &\hspace*{-2ex}\mod m & (\diamond)
\end{array}\right.
\]
where $(\star)$ has solutions \hide{(0 solutions if $m$ is even and $r$ is odd, 1 solution if $m$ is odd and 2
solutions if both $m$ and $r$ are even)}
$q=\frac{m-r}{2}$ if $m$, $r$ have the same parity and
$q=m-\frac{r}{2}$ if $r$ is even.

\medskip

\noindent If $m$ and $r$ are both odd, then $q=\frac{m-r}{2}$ is the only solution of $(\star)$. In this case $(\diamond)$ is equivalent to
$(m-r)(s-1)\equiv 0\mod 2m$; but since $(m-r)(s-1)\equiv 0\mod m$ by \eqref{metacond}, and $m$ is odd, then $(\diamond)$ is
equivalent to $(m-r)(s-1)\equiv 0\mod 2$, which is true. Hence $(\diamond)$ is superfluous in this case.

\smallskip

\noindent If $m$ is odd and $r$ is even, then $q=m-\frac{r}{2}$ is the only solution of $(\star)$. With a similar
argument we see that $(\diamond)$ is superfluous in this case too. \oop{In fact, $(\diamond)$ is equivalent to
$(2m-r)(s-1)\equiv 0 \mod 2m$. By \eqref{metacond}, $(2m-r)(s-1)\equiv 0\mod m$ and since $m$ is odd $(\diamond)$ is
equivalent to $(2m-r)(s-1)\equiv 0\mod 2$ which is true. Hence $(\diamond)$ is superfluous.}

\smallskip

\noindent If $m$ and $r$ are both even, then $(\star)$ has two solutions $q_1=\frac{m-r}{2}$ and
$q_2=m-\frac{r}{2}$. In this case $s$ is odd and $(\diamond)$ is equivalent to $\frac{r}{2}(s-1)\equiv 0\mod m$.\oop{ In fact,
$q_2(s-1)\equiv 0\mod m\Leftrightarrow\frac{r}{2}(s-1)\equiv 0\mod m$ and $q_1(s-1)\equiv 0\mod m\Leftrightarrow\frac{m-r}{2}(s-1)\equiv 0\mod m
\Leftrightarrow (m-r)(s-1)\equiv 0\mod 2m$. Since $s-1$ is even, $m(s-1)\equiv 0\mod 2m$ and so $(m-r)(s-1)\equiv 0\mod 2m
\Leftrightarrow r(s-1)\equiv 0\mod 2m\Leftrightarrow\frac{r}{2}(s-1)\equiv 0\mod m$.} Now the congruence $(\ast)$
says that $\pmeta{m,\frac{n}{2},\frac{r}{2},s}$ is a metacyclic group.

\smallskip

\noindent If $m$ is even and $r$ is odd, then $(\star)$ has no solution.

\medskip

\noindent In the following table we list the central involutions of $G$ according to the parity of the parameters
$m,n,r$ and $s$.

\begin{center}\addtocounter{tab}{1}
\begin{tabular}{r|cccc|l}
\# & $m$ & $n$ & $r$ & $s$ & central involutions\\
\hline
\textbf{1} & odd & odd & odd & odd & ---\\
\textbf{2} & odd & odd & odd & even & ---\\
\textbf{3} & odd & odd & even & odd & ---\\
\textbf{4} & odd & odd & even & even & ---\\
\textbf{5} & odd & even & odd & odd & $y^{\frac{n}{2}}x^{\frac{m-r}{2}}\quad$ if ($\ast$)\\
\textbf{6} & odd & even & odd & even & $y^{\frac{n}{2}}x^{\frac{m-r}{2}}\quad$ if ($\ast$)\\
\textbf{7} & odd & even & even & odd & $y^{\frac{n}{2}}x^{-\frac{r}{2}}$\hspace{3.5ex} if ($\ast$)\\
\textbf{8} & odd & even & even & even & $y^{\frac{n}{2}}x^{-\frac{r}{2}}$\hspace{3.5ex} if ($\ast$)\\
\textbf{9} & even & odd & odd & odd & $x^{m/2}$\\
\textbf{10} & even & odd & even & odd & $x^{m/2}$\\
\textbf{11} & even & even & odd & odd & $x^{m/2}$\\
\textbf{12} & even & even & even & odd & $x^{m/2}$ and\quad $y^{\frac{n}{2}}x^{\frac{m-r}{2}}$,
$y^{\frac{n}{2}}x^{-\frac{r}{2}}$ if ($\ast$) and $\frac{r}{2}(s-1)\equiv 0\mod m$\\
\hline
\end{tabular}\\[5pt]
Table \thetab: The central involutions of a metacyclic group.
\end{center}

\subsection{Presentation pairs}\label{PresPairs}

Let $a=y^px^q,b=y^ux^v\in G=\pmeta{m,n,r,s}$. From \eqref{perm(x,y)} and \eqref{powers}, the pair $(a,b)$ is a
presentation pair for $G$ if and only if $\langle a,b\rangle=G$ and
\[
\begin{array}{l}
pm\equiv 0\mod n\,,\\[2pt]
pr\equiv 0\mod n\,,\\[2pt]
p(s-1)\equiv 0\mod n\,,
\end{array}
\quad
\begin{array}{l}
q\sigma(s^p,m)\equiv 0\mod m\,,\\[2pt]
ru+v\sigma(s^u,n)-q\sigma(s^p,r)\equiv r\frac{pr}{n}\mod m\,,\\[2pt]
v(1-s^p)+q\left( s^u-\sigma(s^p,s)\right)\equiv\frac{pr}{n}(s-1)\mod m\,.
\end{array}
\]
By Theorem \ref{SamePSet} and Corollary \ref{Aut=SP} these equations determine the automorphism group of $G$.

\subsection{Double coverings of cyclic groups}

Every double covering of the cyclic group $C_m=\pmeta{m,1,m,1}$ with presentation $P=\langle x\mid x^m=1\rangle$ of simple
type has presentation
\[
P_1=\langle x,i\mid x^m\rangle\quad\text{or}\quad P_i=\langle x\mid x^m=i\,\rangle\sim\langle x\mid x^{2m}\rangle\,.
\]

\noindent\textbf{Case \textbf{\emph{m}} even.} In this case $[P_1]$ and $[P_i]$ are distinct presentation classes of two non-isomorphic
double coverings of $C_m$, namely $C_m\times C_2$ and $C_{2m}$ respectively. The first double covering is not strong since $x^{m/2}$
is a central involution in $C_m$ (Corollary \ref{StrongGxC2}). On the other hand $C_{2m}$ is a 2-strong double covering of $C_m$
since $[P_i]$ is the only
presentation class of $C_{2m}$ and $C_{2m}$ has only one central involution (Theorem \ref{strong} and Corollary \ref{q-simpletype}).

\hide{
As $x^{m/2}$ is a central involution in
$C_m$ the double covering $C_m\times C_2$ is not strong (Corollary \ref{StrongGxC2}). Since $[P_i]$ is the only
presentation class of $C_{2m}$, which has a unique central involution, $C_{2m}$ is a 2-strong double covering of $C_m$
(Theorems \ref{UniquePClass}, \ref{strong} and Corollary \ref{q-simpletype}).}

\smallskip

\noindent\textbf{Case \textbf{\emph{m}} odd.} In this case there is only one presentation class, namely $[P_1]=\{P_1,P_i\}$, and
therefore we have only one double covering $C_m\times C_2\cong C_{2m}$. Since this has only one central involution, the double covering is 1-strong
(Theorem \ref{strong} and Corollary \ref{q-simpletype}).
\[\addtocounter{diag}{1}
\xymatrix @R=36pt @C=6pt{
C_m\times C_2 \ar@{--}[dr] & & C_{2m} \ar@{-}[dl]|*+{\txt{2-}} & & C_m\times C_2\cong C_{2m} \ar@{-}[d]|*+{\txt{1-}}\\
& {\begin{array}{c} C_m\\ \text{($m$ even)} \end{array}} & & & {\begin{array}{c} C_m\\ \text{($m$ odd)} \end{array}}
}
\]
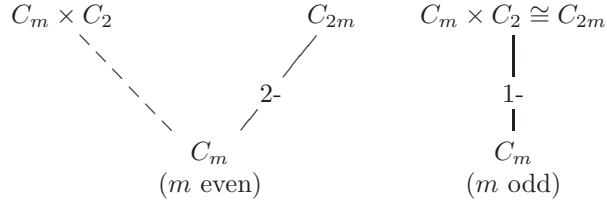
\centerline{Diagram \thediag: The double coverings of cyclic groups.}
%

\subsection{Double coverings of dihedral groups}

Consider the dihedral group $D_m=\pmeta{m,2,m,m-1}$ with presentation $\,P=\langle x,y\mid x^2,y^2,(xy)^m\rangle$\,, where
$m>1$. The transposition of the two generators $x,y$ of $P$ gives rise to an automorphism of
$D_m$.\oop{ In other words, if $(a,b)$ is a presentation pair of $P$ then so does $(b,a)$.} This implies that
$P_{(1,i,\varepsilon)}\sim P_{(i,1,\varepsilon)}$ though $P_{(1,i,\varepsilon)}\not\approx P_{(i,1,\varepsilon)}$, where
$\,\varepsilon=1,i$.

\medskip

\noindent\textbf{Case \textbf{\emph{m}} even.} In this case $P$ is a presentation of simple type of degree $d=0$. It lifts to 8
presentation classes, each containing a single presentation. As seen above, $P_{(1,i,\varepsilon)}\sim
P_{(i,1,\varepsilon)}$ and hence we have at most 6 potentially distinct double coverings of $D_m$.
\begin{itemize}
\item[D1)] $P_1$ is a presentation of $D_m\times C_2$.
\item[D2)] $P_{(1,1,i)}\sim\langle x,y\mid x^2,y^2,(xy)^{2m}\rangle$ which is a presentation of $D_{2m}$.
\item[D3)] $P_{(1,i,1)}\sim P_{(i,1,1)}\sim\left\langle y,z\mid y^4,z^2,(yz)^m,[y^2,z]\right\rangle$. Adding
$x=(yz)^2\,\Leftrightarrow\,y^z=y^{-1}x$ to the presentation we get
\[\left\langle y,z\mid y^4,z^2,(yz)^m,[y^2,z]\right\rangle\sim\left\langle x,y,z\mid x^{\frac{m}{2}},y^4,xx^y,z^2,xx^z,x^{-1}yy^z\right\rangle\cong
\left( \langle x\rangle\rtimes \langle y\rangle\right)\rtimes \langle z\rangle\,.\]
Note that $x^{-1}=(zy^{-1})^2=(zy)^2=x^z=x^y$.\oop{ Note also that $[y^2,z]=y^{-2}(z^{-1}yz)^2=y^{-3}xy^{-1}x=x^{y^{-1}}x$ and that
$[y^2,z]=1\Leftrightarrow x^{y^{-1}}x=1\Leftrightarrow xx^y=1$.} Hence $P_{(1,i,1)}$ is a presentation of
$\left( C_{\frac{m}{2}}\rtimes C_4\right)\rtimes C_2$.
\item[D4)] $P_{(1,i,i)}\sim P_{(i,1,i)}\sim\langle y,z\mid (yz)^my^2,z^2,(yz)^{2m},[y^2,z]\rangle$. Adding
$x=(yz)^2\Leftrightarrow y^z=y^{-1}x$ to the presentation (note that $x^{-1}=(zy^{-1})^2=(zy)^2=x^z=x^y$) we get
\[
P_{(i,1,i)}\sim\langle x,y,z\mid x^{\frac{m}{2}}y^2,x^m,xx^y,z^2,xx^z,x^{-1}yy^z\rangle\,.
\]
This is a presentation of $Q_m\rtimes C_2$, where $Q_{2n}$ is the generalized
quaternion group of order $4n$ with presentation
\[
\langle a,b\mid a^nb^2,a^{2n},aa^b\rangle\sim\langle a,b\mid a^nb^2,aa^b\rangle\,.
\]
\item[D5)] $P_{(i,i,1)}\sim\left\langle x,y\mid x^4,x^2y^2,(xy)^m,[x^2,y]\right\rangle$. Changing generators
($a=xy$ and $b=x$) we get $ P_{(i,i,1)}\sim\left\langle a,b\mid a^m,b^4,aa^b\right\rangle$, which is a
presentation of $C_m\rtimes C_4$.
\item[D6)] $P_{(i,i,i)}\sim\left\langle x,y\mid x^2y^2,(xy)^my^2,y^4\right\rangle$\,. Changing generators ($a=xy$ and $b=y$) we get
$P_{(i,i,i)}\sim\langle a,b\mid a^mb^2,b^4,aa^b\rangle\sim\langle a,b\mid a^mb^2,aa^b\rangle$ which is a presentation
of $Q_{2m}$.
\end{itemize}
\textbf{Case \textbf{\emph{m}} odd.} In this case we have $ P_{(j_1,j_2,1)}\sim P_{(j_2,j_1,1)}\approx P_{(j_2,j_1,i)}\sim
P_{(j_1,j_2,i)}$, for every $j_1,j_2\in\{1,i\}$. Moreover, for every $(j_1,j_2,j_3)\in\{1,i\}^3$ we have
\[
\begin{array}{ll}
j_3 & =(xy)^m=(j_1x^{-1}j_2y^{-1})^m=j_1j_2\,(yx)^{-m}=j_1j_2\,\left((xy)^{y^{-1}}\right)^{-m}\\
& =j_1j_2\,\left((xy)^{-m}\right)^{y^{-1}}=j_1j_2\,j_3^{y^{-1}}=j_1j_2j_3\,,
\end{array}
\]
that is $j_1j_2=1$. If $j_1\neq j_2$ then $i=1$ and $P_{(j_1,j_2,j_3)}$ is a presentation of $D_m$, hence not a double covering.
Thus $j_1=j_2$ and we have two cases. Namely,

\begin{itemize}
\item[D1)] $P_{(1,1,1)}\approx P_{(1,1,i)}$ are presentations of $D_m\times C_2\cong D_{2m}$.
\item[D2)] $P_{(i,i,1)}\approx P_{(i,i,i)}$ are presentations of $Q_{2m}$. This can be easily seen by
changing generators $a=xy$ and $b=y$.
\end{itemize}

If $m=2$ then $\ D_{2m}\cong\left( C_{\frac{m}{2}}\rtimes C_4\right)\rtimes C_2\cong Q_m\rtimes
C_2\cong C_m\rtimes C_4$ in which case we have 3 non-isomorphic double coverings of $D_2=V_4$, namely $D_2\times C_2$, $D_4$ and $Q_4$.
While the first two are not strong, $Q_4$ is a 4-strong double covering.

\smallskip

For $m>2$ we have 6 or 2 non-isomorphic double coverings of $D_m$, according as $m$ is even or odd.
As both $D_{2m}=\pmeta{2m,2,2m,2m-1}$ and $Q_{2m}=\pmeta{2m,2,m,2m-1}$ have only one central involution (Table 1, lines
\textbf{11} and \textbf{12}) and only one presentation class, they are
$q$-strong double coverings of $D_m$, where $q=2$ if $m$ is odd and $q=4$ if $m$ is even (Theorem \ref{strong}
and Corollary \ref{q=SizeKrho}). If $m$ is odd, these are the only double coverings of $D_m$.

\smallskip

Let $m$ be even. Then $D_m$ has a central involution and therefore $D_m\times C_2$ is not a strong double covering of $D_m$
(Corollary \ref{StrongGxC2}). The double coverings $\left( C_{\frac{m}{2}}\rtimes C_4\right)\rtimes C_2$ and
$Q_m\rtimes C_2$ are not strong since both have more then one presentation class (Theorem \ref{strong}).
The double covering $\hat G=C_m\rtimes C_4$, having presentation
\[
\langle x,y\mid x^m,y^4,xx^y\rangle\sim\langle x,y\mid x^m=1,y^4=x^m,x^y=x^{-1}\rangle=\pmeta{m,4,m,m-1}\,,
\]
is a metacyclic group with central involutions $x^{m/2}$, $y^2x^{m/2}$ and $y^2$ (Table 1, \#\textbf{12}). Now
$[P_{(i,i,1)}]$ is the only presentation class of $\hat G$. Then $\hat G$ is strong (hence 4-strong by Corollary
\ref{q-simpletype})\hide{double covering of $D_m$} if and only if $i=y^2$ is fixed by every automorphism of $\hat G$
(Theorem \ref{strong}). If $\psi\in Aut(\hat G)$ maps $(x,y)$ to $(a,b)=(y^px^q,y^ux^v)$ then $(a,b)$ is another
presentation pair for $\pmeta{m,4,m,m-1}$. By \S\ref{PresPairs},
\[
p\,m\equiv 0\mod 4\qquad\text{and}\qquad p\,(m-2)\equiv 0\mod 4
\]
which is equivalent to $p$ even (0 or 2). Then $\langle a,b\rangle =\hat G$ implies $u$ odd (1 or 3). Thus
\[
i\psi=y^2\psi=b^2=(y^ux^v)^2=y^{2u}x^{v(1+(-1)^u)}=y^{2u}=y^2=i\,,
\]
proving that $\hat G=C_m\rtimes C_4$ is a strong double covering of $D_m$.

\[\addtocounter{diag}{1}
\xymatrix @C=0pt @R=36pt{ D_m\!\!\times\! C_2 \ar@{--}[drrr] & D_{2m} \ar@{-}[drr]|*+{\txt{4-}} &
(C_{\frac{m}{2}}\!\!\rtimes\! C_4)\!\!\rtimes\! C_2 \ar@{--}[dr] & Q_m\!\!\rtimes\! C_2\ar@{--}[d] & C_m\!\!\rtimes\!
C_4 \ar@{-}[dl]|*+{\txt{4-}} & Q_{2m} \ar@{-}[dll]|*+{\txt{4-}} &
\quad D_m\!\!\times\! C_2\cong D_{2m}\ar@{-}[dr]|*+{\txt{2-}} & & Q_{2m} \ar@{-}[dl]|*+{\txt{2-}}\\
& & & {\begin{array}{c} D_m\\ \text{($m>2$ even)}\end{array}} & & & & {\begin{array}{c} D_m\\ \text{($m$
odd)}\end{array}} }
\]
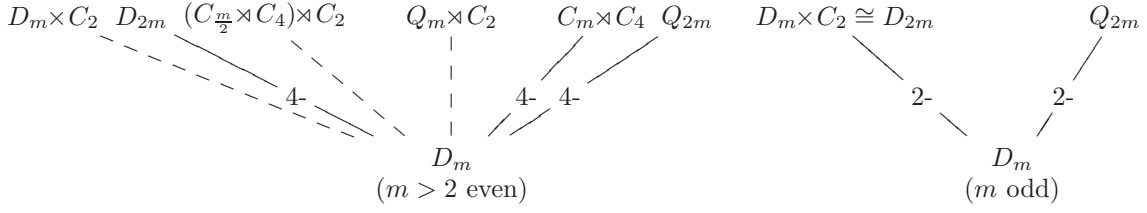
\centerline{Diagram \thediag: The double coverings of the dihedral group $D_m$.}
%

\begin{rem}\label{FixAllInv}
We have seen above that the central involution $i=y^2$ of $C_m\rtimes C_4$ is fixed by any automorphism. However,
if $m/2$ is even, using \eqref{powers} and taking in account that for this case $\sigma(s^p,k)\equiv
k\mod m$, we also have
\[
x^{\frac{m}{2}}\psi=(x\psi)^{\frac{m}{2}}=a^{\frac{m}{2}}=(y^px^q)^{\frac{m}{2}}=y^{p{\frac{m}{2}}}x^{q{\frac{m}{2}}}=
x^{q{\frac{m}{2}}}=x^{\frac{m}{2}}\,.
\]
The last equality follows from the fact that the order of $x^{\frac{m}{2}}\psi$ is 2. So, for $m\equiv 0\mod 4$ all the
three central involutions of $\pmeta{m,4,m,m-1}$ are fixed by any automorphism of $\pmeta{m,4,m,m-1}$.
\end{rem}

\subsection{Double coverings of dicyclic groups}

A dicyclic group, often called ``generalised quaternion group'', is a group $Q_{2m}=\pmeta{2m,2,m,2m-1}$ of order $4m$ with
presentation
\[P=\langle x,y\mid x^my^{-2}, xx^y\rangle\,.\]
This presentation has deficiency zero and is of simple type of degree 0 or 1 according as $m$ is even or odd.

\smallskip

\noindent\textbf{Case \textbf{\emph{m}} even.} In this case $P$ lifts to 4 singular presentation classes $[P_1]$,
$[P_{(i,1)}]$, $[P_{(1,i)}]$ and $[P_{(i,i)}]$.
\begin{itemize}
\item[Q1)] $P_1$ is a presentation of $Q_{2m}\times C_2$.
\item[Q2)] $P_{(i,1)}\sim\langle x,y\mid x^{2m},y^4,xx^y\rangle\,$ which is a presentation of $C_{2m}\rtimes C_4$.
\item[Q3)] $P_{(1,i)}\sim P_{(i,i)}$, by changing generators $a=x$, $b=xy$ and $i=i$,\oop{ (in fact, $P_{(1,i)}\sim
\langle a,b,i\mid a^m(a^{-1}b)^{-2},aa^bi,i^2,[i,a],[i,b]\rangle\sim P_{(i,i)}\,$)} which is a presentation of the
binary group $\widetilde{Q}_{2m}=P_{(i,i)}$.
\end{itemize}

\noindent\textbf{Case \textbf{\emph{m}} odd.} In this case $P$ lifts to 2 presentation classes.
\begin{itemize}
\item[Q1)] $P_1\approx P_{(i,1)}$ are presentations of $Q_{2m}\times C_2$.
\item[Q2)] $P_{(1,i)}\approx P_{(i,i)}$ are presentation of the binary group $\widetilde{Q}_{2m}=P_{(i,i)}$. Since
\[
\begin{array}{ll}
P_{(1,i)} & =\langle x,y\mid x^m=y^2,xx^y=i\rangle\\
& \sim\langle x,y\mid x^m=y^2,xx^y=i,x^{2m}=i\rangle\\
& \sim\langle x,y\mid x^{4m}=1, y^2=x^m, x^y=x^{2m-1}\rangle
\end{array}
\]
we see that the binary $\widetilde{Q}_{2m}$ is a metacyclic group $\pmeta{4m,2,m,2m-1}$.
\end{itemize}

\noindent Since $Q_{2m}$ has a central involution (Table 1, \#\textbf{11} and \textbf{12}), $Q_{2m}\times C_2$ is
not strong (Corollary \ref{StrongGxC2}). If $m$ is even then $\widetilde{Q}_{2m}$ is not strong as well since it has
two presentation classes (Theorem \ref{strong}). If $m$ is odd then $\widetilde{Q}_{2m}$ is a metacyclic group
$\pmeta{4m,2,m,2m-1}$ which has a unique central involution (Table 1, \#\textbf{11}) and is therefore 2-strong (Theorem
\ref{strong} and Corollary \ref{q-simpletype}). Finally, if $m$ is even, $C_{2m}\rtimes C_4\not\cong\widetilde{Q}_{2m}$
is a 4-strong double covering of $Q_{2m}$ (Remark \ref{FixAllInv}, Theorem \ref{strong} and Corollary
\ref{q-simpletype}).

\[\addtocounter{diag}{1}
\xymatrix @C=6pt @R=36pt{ Q_{2m}\!\!\times\! C_2 \ar@{--}[dr] & C_{2m}\!\!\rtimes\! C_4 \ar@{-}[d]|*+{\txt{4-}} &
\tilde{Q}_{2m} \ar@{--}[dl] &
& Q_{2m}\!\!\times\! C_2 \ar@{--}[dr] & & \tilde{Q}_{2m} \ar@{-}[dl]|*+{\txt{2-}}\\
& {\begin{array}{c} Q_{2m}\\ \text{($m$ even)}\end{array}} & & & & {\begin{array}{c} Q_{2m}\\ \text{($m$
odd)}\end{array}} }
\]
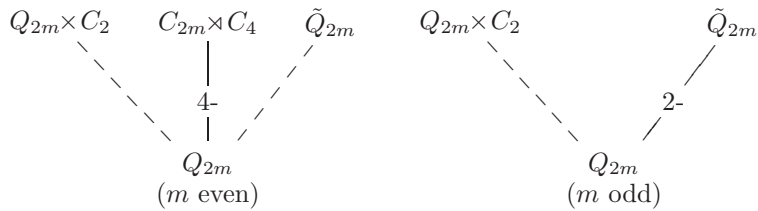
\centerline{Diagram \thediag: The double coverings of the dicyclic groups.}
%

\subsection{Double coverings of metacyclic groups: The general case}

Let $G$ be a metacyclic group with presentation $P=\pmeta{m,n,r,s}$ where the parameters $m,n,r,s\in\mathds{N}$ satisfy
\eqref{metacond} and assume $r,s\leqslant m$. In the following, for $\epsilon\in\{0,1\}$, let $i^\epsilon\in C_2$ have the usual meaning $i^0=1$ and
$i^1=i$. For every $J=(i^\epsilon,i^\tau,i^\nu)\in C_2^3$, where $C_2=\{1,i\}$, $P_J=\langle x,y\mid x^m=i^\epsilon,y^n=x^ri^\tau,x^y=x^si^\nu\rangle$
with $m$, $n$, $r$ and $s$ satisfy \eqref{metacond}. This
is either a presentation of a double covering of $G$ or a presentation of $G$. In the following theorems we give necessary and
sufficient conditions that assures when $P_J$ is a presentation of a double covering of $G$ or not, for every possible
choice of $J\in C_2^3$. If $J=1$ there are no conditions since $P_1$ is a presentation of $G\times C_2$ and therefore a
presentation of a double covering of $G$.

\medskip

\begin{thm}\label{Pixx}
Let $J=(i,i^\epsilon,i^\tau)\in C_2^3$, where $\epsilon,\tau\in\{0,1\}$. The following statements are equivalent:
\begin{itemize}
\item[\emph{(a)}] $P_J$ is a double covering of $\,G$\,.
\item[\emph{(b)}] $P_J$ is a metacyclic group $\,\pmeta{2m,n,r+m\epsilon,s+m\tau}$\,.
\item[\emph{(c)}] $(s+m\tau)^n\equiv 1\hspace{-1ex}\mod 2m\quad\text{and}\quad (r+m\epsilon)(s+m\tau)\equiv r+m\epsilon\!\mod 2m$\,.
\end{itemize}
\end{thm}

\begin{proof}
It is enough to prove that (a) implies (b).
Since $\ P_J\sim\langle x,y\mid x^{2m}=1,y^n=x^{r+m\epsilon}, x^y=x^{s+m\tau},[x^m,y]=1\rangle\,$ $P_J$ is clearly a metacyclic group of
order $2mn$. Hence the relation $[x^m,y]=1$ is superfluous and $P_J=\pmeta{2m,n,r+m\epsilon,s+m\tau}$.
\oop{(a) implies $x^m\neq 1$ i.e. $x^m$ is an involution. From the third relation
$x=x^{y^n}=x^{(s+m\tau)^n}\Rightarrow (s+m\tau)^n\equiv 1\mod 2m\Rightarrow s+m\tau$ is odd. From the third relation again
$(x^m)^y=(x^m)^{s+m\tau}=x^m$\,. Hence $[x^m,y]=1$ is superfluous.}
\hide{then (a) implies that
$x$ has order $2m$ and this implies (b). If (b) is satisfied then $P_J=\pmeta{2m,n,r+m\epsilon,s+m\tau}/\overline{\langle [x^m,y]\rangle}$.
Since in this case
$s+m\tau$ is odd,\, $x^y=x^{s+m\tau}$ implies $[x^m,y]=1$. Hence $P_J$ is the metacyclic group
$\pmeta{2m,n,r+m\epsilon,s+m\tau}$.}
\end{proof}

\begin{cor}
Let $J=(i,1,i^\tau)$ and $K=(i,i,i^\tau)$, where $\tau\in\{0,1\}$. Then $P_J$ is a double covering of $\,G$ if and only
if $P_K$ is a double covering of $\,G$.
\end{cor}

\begin{proof}
This follows from statement (b) of previous theorem. Assuming $(s+m\tau)^n\equiv 1\mod 2m$, which
implies $s+m\tau$ odd, we have that $r(s+m\tau)\equiv r\mod 2m$ is equivalent to $(r+m)(s+m\tau)\equiv r+m\mod 2m$\,.
\end{proof}

\begin{cor}\label{cPixx}
Let $J=(i,i^\epsilon,1)$ and $K=(i,i^\epsilon,i)$ where $\,\epsilon\in\{0,1\}\,$. Then $P_J$ and $P_K$ are both double coverings of $\,G$
if and only if $\,s^n\equiv 1\mod 2m$, $rs\equiv r\mod 2m$ and $m,n,r$ are even.
\end{cor}

\begin{proof}
From previous theorem we have that $P_J$ and $P_K$ are both double coverings of $G$
if and only if
\[
\begin{array}{lll}
\text{(a)} & s^n\equiv 1 & \mod 2m,\\[2pt]
\text{(b)} & (r+m\epsilon)s\equiv r+m\epsilon & \mod 2m,\\[2pt]
\text{(c)} & (s+m)^n\equiv 1 & \mod 2m,\\[2pt]
\text{(d)} & (r+m\epsilon)(s+m)\equiv r+m\epsilon & \mod 2m.
\end{array}
\]
Assuming (a), which implies $s$ odd, then (b) is equivalent to $rs\equiv r\mod 2m$. As (c) implies
$s+m$ odd, it follows from (a) and (c) that $m$ is even. It also comes from (a) and (c) that $n$ is even, since $n>1$ and
\[
(s+m)^n=s^n+mns^{n-1}+m^2\kappa\,,\qquad \text{for some }\kappa\in\mathds{N}\,.\]
On the other hand, (b) and (d) imply $r+m\epsilon$ even which together with $m$ even implies $r$ even.
Reciprocally, (a), (b) and $m,n,r$ even imply (c) and (d).
\end{proof}

\begin{thm}\label{Plxi}
Let $J=(1,i^\epsilon,i)\in C_2^3$ where $\epsilon\in\{0,1\}$. Then $P_J$ is a double covering of $\,G$ if and only if
$\,m$, $n$ and $r$ are even, in which case $P_J$ is a metabelian group of order $2mn$.
\end{thm}

\begin{proof}
($\Rightarrow$): Suppose that $\,P_J=\langle x,y\mid x^m=1, y^n=x^ri^\epsilon, x^y=x^si\rangle$\, is a double covering of $G$. Then $i\neq 1$ and from
condition \eqref{metacond} we get
\[
\begin{array}{l}
1\oop{=(x^m)^y}=(x^y)^m=(x^si)^m=i^m\,,\\
x=x^{y^n}=x^{s^n}i^{\sigma(s,n)}=x\,i^{\sigma(s,n)},\qquad\text{where}\quad\sigma(s,n)=\sum\limits_{k=0}^{n-1}s^k\,,\\
x^r\oop{=(x^r)^y}=(x^y)^r=(x^si)^r=x^{rs}i^r=x^r\,i^r\,.
\end{array}
\]
Then $m$, $\sigma(s,n)$ and $r$ are even. Hence, by the condition \eqref{metacond}, $s$ is odd and so $i^{\sigma(s,n)}=1$ is equivalent to $i^n=1$ which
implies that $n$ is even.

\smallskip

\noindent ($\Leftarrow$): Since $N=\langle x,i\rangle$ is a normal subgroup of $P_J$ and $P_J/N\cong C_n$ we have that $|P_J|=|N|n$\,.
Taking into account that $m$, $n$, $r$ even implies $\sigma(s,n)$ even, since by condition \eqref{metacond} $s$ is odd, by the
Reidemeister-Schreier rewriting process we get that $\langle x,i\mid x^m, i^2, [x,i]\rangle$ is a presentation for $N$\,.
Hence $N\cong C_m\times C_2$\,. This proves that $P_J$ is a metabilian group of order $2mn$ and hence a double covering of $G$.
\end{proof}

\begin{thm}\label{Plil}
$P_{(1,i,1)}$ is the metacyclic group $\pmeta{m,2n,2r,s}$ and hence a double covering of $\,G$.
\end{thm}

\begin{proof}
Taking into account \eqref{metacond} we get
$\,P_{(1,i,1)}\sim\langle x,y\mid x^m=1,y^{2n}=x^{2r}, x^y=x^s\rangle$\,.
As $s^n\equiv 1\mod m$ implies $s^{2n}\equiv 1\mod m$ and as $rs\equiv r\mod m$ implies $2rs\equiv 2r\mod m$,
$P_{(1,i,1)}$ is the metacyclic group $\pmeta{m,2n,2r,s}$ of order $2mn$ and therefore a double covering of $G$.
\end{proof}

\begin{cor}\label{cPlil}
The following statements are equivalent:
\begin{itemize}
\item[\emph{(a)}] For all $J\in C_2^3$, $P_J$ is a presentation of a double covering of $G$\,.
\item[\emph{(b)}] $m,n,r$ are even and $\,s^n\equiv 1\mod 2m$, $rs\equiv r\mod 2m$\,.
\end{itemize}
\end{cor}

\begin{proof}
The statement follows from Corollary \ref{cPixx} and
Theorems \ref{Plxi} and \ref{Plil}.
\oop{Notice that $P_1$ is a presentation of the double covering $G\times C_2$ of $G$.}
\end{proof}

\begin{lem}
If $\,(s+m\tau)^n\equiv 1\mod m$ holds for $\tau=0,1$ and $m$ is even, then $\frac{s^n-1}{m}$ and $\frac{(s+m)^n-1}{m}$ have
the same parity if and only if $\,n$ is even.
\end{lem}

\begin{proof}
Since $(s+m)^n=s^n+mns^{n-1}+\kappa m^2$ for some $\kappa\in\mathds{N}$ we can write
$\frac{(s+m)^n-1}{m}=\frac{s^n-1}{m}+ns^{n-1}+\kappa m$. Now, since $\kappa m$ is even and $ns^{n-1}$ has the same parity as $n$ (as $s$ is odd),
then $ns^{n-1}+\kappa m$ has the same parity as $n$. Hence $\frac{(s+m)^n-1}{m}$ and $\frac{s^n-1}{m}$ have the same
parity if and only if $n$ is even.
\end{proof}

\medskip

\noindent This lemma together with Theorem \ref{Pixx} gives the following corollary.

\begin{cor}\label{Simultan}
If $\,m,n$ are even and $\frac{s^n-1}{m}$ is odd, then for any $J=(i,i^\epsilon,i^\tau)$, where $\epsilon,\tau\in\{0,1\}$,
$P_J$ is not a double covering of $\,G$\,.
\end{cor}

Table below displays the presentation classes $[P_J]$ of possible double coverings of the metacyclic group $G$ with
presentation $P=\pmeta{m,n,r,s}$ according to the parity of the parameters. The number of such presentation classes is an
upper bound of the number of (non-isomorphic) double coverings of $G$ (Theorem \ref{NDC}). On the right column we
display a better upper bound $\delta$ based on Theorems \ref{Pixx}, \ref{Plxi}, \ref{Plil} and Corollary \ref{cPixx}.
The erased presentation classes are, according to Theorem \ref{Plxi}, presentation classes of $G$ and therefore not
presentation classes of a double covering of $G$.
\def\strike{\raisebox{4pt}{\rule{35pt}{0.5pt}}\hspace{-34pt}}
\begin{center}\addtocounter{tab}{1}
\begin{math}
\begin{array}{l|cccc|l|l}
\# & m & n & r & s & [P_J] & \delta \\ \hline
   & & & & & & \\[-8pt]
\mathbf{1} & \text{odd} & \text{odd} & \text{odd} & \text{odd} & [P_1]=[P_{(1,i,1)}],\ \strike [P_{(1,1,i)}] & 1 \\[2pt]
\mathbf{2} & \text{odd} & \text{odd} & \text{odd} & \text{even} & [P_1]=[P_{(1,i,1)}],\ \strike [P_{(1,i,i)}] & 1 \\[2pt]
\mathbf{3} & \text{odd} & \text{odd} & \text{even} & \text{odd} & [P_1]=[P_{(1,i,1)}],\ \strike [P_{(1,1,i)}] & 1 \\[2pt]
\mathbf{4} & \text{odd} & \text{odd} & \text{even} & \text{even} & [P_1]=[P_{(1,i,1)}],\ \strike [P_{(1,i,i)}] & 1 \\[2pt]
\mathbf{5} & \text{odd} & \text{even} & \text{odd} & \text{odd} & [P_1],\ [P_{(1,i,1)}],\ \strike [P_{(1,i,i)}],\ \strike [P_{(1,1,i)}] & 2\\[2pt]
\mathbf{6} & \text{odd} & \text{even} & \text{odd} & \text{even} & [P_1],\ [P_{(1,i,1)}],\ \strike [P_{(1,i,i)}],\ \strike [P_{(1,1,i)}] & 2 \\[2pt]
\mathbf{7} & \text{odd} & \text{even} & \text{even} & \text{odd} & [P_1],\ [P_{(1,i,1)}],\ \strike [P_{(1,1,i)}],\ \strike [P_{(1,i,i)}] & 2 \\[2pt]
\mathbf{8} & \text{odd} & \text{even} & \text{even} & \text{even} & [P_1],\ [P_{(1,i,1)}],\ \strike [P_{(1,1,i)}],\ \strike [P_{(1,i,i)}] & 2 \\[2pt]
\mathbf{9} & \text{even} & \text{odd} & \text{odd} & \text{odd} & [P_1]=[P_{(1,i,1)}],\ [P_{(i,1,1)}],\ [P_{(i,1,i)}],\ \strike [P_{(1,1,i)}] & 2 \\[2pt]
\mathbf{10} & \text{even} & \text{odd} & \text{even} & \text{odd} & [P_1]=[P_{(1,i,1)}],\ [P_{(i,1,1)}],\ [P_{(i,1,i)}],\ \strike [P_{(1,1,i)}] & 2 \\[2pt]
\mathbf{11} & \text{even} & \text{even} & \text{odd} & \text{odd} & [P_1]=[P_{(1,i,1)}],\ [P_{(i,1,1)}],\ [P_{(i,1,i)}],\ \strike [P_{(1,1,i)}] & 2 \\[2pt]
\mathbf{12} & \text{even} & \text{even} & \text{even} & \text{odd} & [P_1],[P_{J_2}],\dots, [P_{J_8}],\ \ J_i\in C_2^3\setminus\{1\}& 8 \\ \hline
\end{array}
\end{math}\\[10pt]
Table \thetab: The double coverings of metacyclic groups.
\end{center}

\begin{rem}
In \#\textbf{1-4} there is exactly one double covering, $G\times C_2$.
In \#\textbf{5-8} there are exactly two double coverings, namely $G\times C_2$ and the metacyclic group
$\pmeta{m,2n,2r,s}$ (Theorem \ref{Plil}) which are non-isomorphic according to Theorem \ref{UniquePClass}. Hence the upper bound $\delta=2$ is always
attained.\\
In \#\textbf{9-11} Corollary \ref{cPixx} justifies the upper bound $\delta=2$ which is also always attained.
In fact, $m$ even implies
\begin{equation}\label{sLines9-11}
\frac{(s+m)^n-1}{m}=\frac{s^n-1}{m}+ns^{n-1}+2\kappa\,,
\end{equation}
for some non-negative integer $\kappa$. Multiplying \eqref{sLines9-11} by $r$ we get
\begin{equation}\label{nLines9-11}
r\frac{(s+m)^n-1}{m}\equiv\frac{r(s+m-1)}{m}n\mod 2\,,
\end{equation}
since $r\frac{s^n-1}{m}=\frac{r(s-1)}{m}\sum\limits_{k=0}^{n-1}s^k\equiv\frac{r(s-1)}{m}n\mod 2$. Now,
\eqref{sLines9-11} and \eqref{nLines9-11} together with the following equality
\begin{equation}\label{rLines9-11}
\frac{r(s+m-1)}{m}=\frac{r(s-1)}{m}+r
\end{equation}
guarantee in \#\textbf{9-11} that $(s+m\tau)^n\equiv 0\mod 2m$ and $r(s+m\tau)\equiv r\mod 2m$\,, for some
$\tau\in\{0,1\}$. According to Theorem \ref{Pixx} this shows that $\delta=2$ is always attained.

\medskip

We have just seen that the upper bound $\delta$ in \#\textbf{1-11} is the exact number of
non-isomorphic double coverings of $G$. The same can not be said in \#\textbf{12}. For example the
dihedral groups $D_m$ and the dicyclic groups $Q_{2m}$ for $m$ even falls in this line and the
number of their double coverings is 6 and 3 respectively.
\end{rem}

The direct product $G\times C_2$ is the only double covering of $G$ in \#\textbf{1-4} and since it
has presentation $P_{(1,i,1)}$ it is a metacyclic group $\pmeta{m,2n,2r,s}$ (Theorem \ref{Plil}).
As $G$ has no central involutions (Table 1, \#\textbf{1-4}), $G\times C_2$ is 1-strong
(Corollaries \ref{StrongGxC2} and \ref{q=SizeKrho}).

In \#\textbf{5-8}, $G\times C_2$ has presentation $P_1\approx P_{(i,i^\epsilon,i^\tau)}$ for some
$\epsilon,\tau\in\{0,1\}$ with $r+m\epsilon$ even, hence $G\times C_2$ is the metacyclic group $\pmeta{2m,n,r+m\epsilon,s+m\tau}$
(Theorem \ref{Pixx}). Now $G$ has a central involution if and only if $s^{\frac{n}{2}}\equiv 1\mod m$
(Table 1, \#\textbf{5-8}), hence $G\times C_2$ is a strong double covering of $G$ if and only if
$s^{\frac{n}{2}}\not\equiv 1\mod m$ (Corollary \ref{StrongGxC2}), being 2-strong
in such case (Corollary \ref{q=SizeKrho}). The other double covering of $G$ has presentation $P_{(1,i,1)}$, which is a metacyclic group
$\pmeta{m,2n,2r,s}$ (Theorem \ref{Plil}) with one central involution only (Table 1, \#\textbf{7} and \#\textbf{8})
and therefore a 2-strong double covering (Theorem \ref{UniquePClass}, Theorem \ref{strong} and Corollary \ref{q=SizeKrho}).

In \#\textbf{9-11}, $G\times C_2$ is the metacyclic group $\pmeta{m,2n,2r,s}$ with presentation
$P_{(1,i,1)}$. As $G$ has a central involution (Table 1, \#\textbf{9-11}), \mbox{$G\times C_2$} is
not strong (Corollary \ref{StrongGxC2}). There is one more double covering $\hat G$ of $G$ other
than $G\times C_2$. This has presentation either $P_{(i,1,1)}$ or $P_{(i,1,i)}$ (Corollary
\ref{cPixx}), being $\hat G$ a metacyclic group, either $\pmeta{2m,n,r,s}$ or $\pmeta{2m,n,r,s+m}$
(Theorem \ref{Pixx}). In either cases, $\hat G$ has only one central involution (Table 1,
\#\textbf{9-11}) which makes it a 2-strong double covering of $G$ (Theorem \ref{UniquePClass},
Theorem \ref{strong} and Corollary \ref{q=SizeKrho}).

In \#\textbf{12}, $G$ has one or three central involutions (Table 1, \#\textbf{12}) and so $G\times
C_2$ is not strong. The others 7 (at most) double coverings $\hat G\not\cong G\times C_2$ have 1 or more central involutions
depending on the parameters. As remarked below, they may have one or
more presentation classes. In this case the strongness of $\hat G$ depends sharply on the choice of the
parameters.

\begin{rem}
In \#\textbf{12} a little more can be said. According to Corollary \ref{Simultan}, if $\frac{s^n-1}{m}$ is odd then
$\delta=4$. The value of $\delta$ cannot be always reduced. Using GAP
\cite{GAP} we see that for $(m,n,r,s)=(10,4,10,3)$ and $(10,4,10,7)$ the eight resulting groups are all non-isomorphic and so the upper bound
$\delta=8$ is reached in these two cases. Observe yet that $\pmeta{m,n,r,s}$ has, according to \cite{Wam}, zero deficiency if and only if
\[gcd\left( m,r,s-1,\frac{s^n-1}{m},\frac{r(s-1)}{m},\frac{s^n-1}{s-1}\right) =1\,.\] In the particular case when $r=\frac{m}{(m,s-1)}$,
the metacyclic group $\pmeta{m,n,r,s}$ has the following presentation of zero deficiency (see Johnson \cite{DJ}, pg 91)
$$\langle x,y\mid y^n=x^r,[y,x^{-t}]=x^{(m,s-1)}\rangle\,,$$ for a certain integer $t$. In this case Table \thetab\ reduces to
\[\begin{array}{c|ccc|c}
\# & n & r & (m,s-1) & \delta \\
\hline
           & & & & \\[-8pt]
\mathbf{1} & \text{odd} & \text{odd} & \text{odd} & 1 \\
\mathbf{2} & \text{odd} & \text{odd} & \text{even} & 2 \\
\mathbf{3} & \text{odd} & \text{even} & \text{odd} & 1 \\
\mathbf{4} & \text{odd} & \text{even} & \text{even} & 2 \\
\mathbf{5} & \text{even} & \text{odd} & \text{odd} & 2 \\
\mathbf{6} & \text{even} & \text{odd} & \text{even} & 2 \\
\mathbf{7} & \text{even} & \text{even} & \text{odd} & 2 \\
\mathbf{8} & \text{even} & \text{even} & \text{even} & 4 \\
\hline
\end{array}\]
\begin{center}\addtocounter{tab}{1}
Table \thetab: A particular case.
\end{center}
As we can see, \#\textbf{12} of Table 2, which corresponds to \#\textbf{8} in Table 3,
has the upper bound $\delta$ dropped to 4.
\end{rem}

\section{Double coverings of rotary platonic groups}\label{DCPG}

By a \emph{rotary platonic group} we mean the rotation group of some platonic solid. Thus the rotary platonic groups are
$A_4$, $S_4$ and $A_5$. For $n=3,4,5$,
\[
P=\langle x,y\mid x^3,y^n,(xy)^2\rangle
\]
is a presentation of simple type of $A_4$, $S_4$ and $A_5$, respectively.

\begin{enumerate}
\item Double coverings of $G=A_4$, $A_5$ ($n=3,5$).

\smallskip

Being $n$ odd there are 2 non-isomorphic double
coverings (Corollary \ref{NDC} and Theorems \ref{UniquePClass} and \ref{SimpleType}), the direct product $G\times C_2$ with presentations $P_{(j_1,j_2,1)}$ and the binary group $\tilde G$ with
presentations $P_{(j_1,j_2,i)}$, where $j_1,j_2\in\{1,i\}$.
\item Double coverings of $S_4$ ($n=4$).

\smallskip

In this case we have $4={8/2}$ presentation classes:

1) $[P_{(1,1,1)}]=[P_{(i,1,1)}]$, giving presentations of $S_4\times C_2$.

2) $[P_{(1,i,1)}]=[P_{(i,i,1)}]$, where $P_{(i,i,1)}\sim\left\langle x,y\mid x^6=(xy)^2=1,x^3=y^4\right\rangle$ giving
presentations of $GL(2,3)$.

3) $[P_{(1,1,i)}]=[P_{(i,1,i)}]$, giving presentations of a group $B$ described below.

4) $[P_{(1,i,i)}]=[P_{(i,i,i)}]$, giving presentations of the binary $\tilde S_4$.

\bigskip

\noindent\textbf{The group B}

\smallskip

In contrast to $GL(2,3)$ and the binary octahedral group $\tilde S_4$ that are double coverings $\hat G$ of $S_4$ with $C_2$
in the derived group $\hat G'$, the group $B$ has $C_2$ not in $B'$, and while for $GL(2,3)$ and $\tilde S_4$ we have
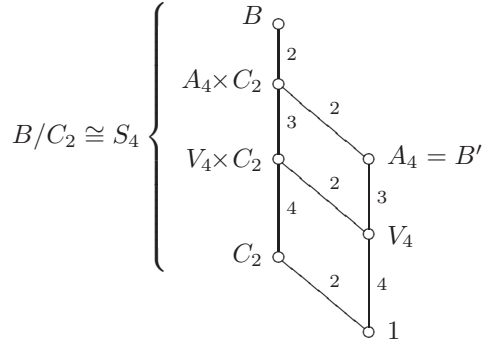
$\hat G/\hat G'\cong C_2$ and for $S_4\times C_2$ we have $\hat G/\hat G'\cong V_4$, for $B$ we have $B/B'\cong C_4$.
The normal subgroups of $B=\langle x,y\mid x^6,y^4,(xy)^2x^3,[x^3,y]\rangle$ are:
$A_4\times C_2=\langle x,y^{-1}xy^{-1}\rangle$, $A_4=B'=\langle x^4,y^{-1}xy^{-1}\rangle$, $V_4\times C_2=\langle
x^3y^2,x^2y^2,(xy)^2\rangle$, $V_4=\langle x^3y^2,x^2y^2x\rangle=\langle x^3y^2\rangle^B$ and $C_2=\langle
(xy)^2\rangle$.
\[\addtocounter{diag}{1}
B/C_2\cong S_4\left\{\begin{array}{r} B\\[12pt] A_4\!\!\times\! C_2\\[18pt] V_4\!\!\times\! C_2\\[24pt] C_2\end{array}\right.\hspace{-4pt} \xy
\POS(0,15)*\cir<2pt>{}="l1"\POS(0,7)*\cir<2pt>{}="l2"\POS(0,-3)*\cir<2pt>{}="l3"\POS(0,-16)*\cir<2pt>{}="l4"
\POS(12,-3)*\cir<2pt>{}="r1"\POS(12,-13)*\cir<2pt>{}="r2"\POS(12,-26)*\cir<2pt>{}="r3"
\POS"l1"\ar@{-}^2"l2"\POS"l2"\ar@{-}^2"r1"\POS"l2"\ar@{-}^3"l3"\POS"r1"\ar@{-}^3"r2"\POS"l3"\ar@{-}^2"r2"\POS"l4"\ar@{-}^2"r3"\POS"l3"\ar@{-}^4"l4"
\POS"r2"\ar@{-}^4"r3"\endxy\hspace{-8pt}
\begin{array}{l}
\\[70pt] A_4=B'\\[18pt] V_4 \\[24pt] 1
\end{array}
\]
\centerline{Diagram \thediag: The normal subgroups of $B$.}
\end{enumerate}
Since $A_4$, $A_5$ and $S_4$ are centerless, every double covering has exactly one central involution which is
obviously fixed by every automorphism. By Theorem \ref{strong}, every double covering of $A_4$, $A_5$ and $S_4$ is
$q$-strong, where by Corollary \ref{q-simpletype}, $q=1$ for double coverings of $A_4$ and $A_5$, and $q=2$ for
double coverings of $S_4$.

\[\addtocounter{diag}{1}
\xymatrix @R=36pt @C=6pt{
A_4\!\!\times\! C_2 \ar@{-}[dr]|*+{\txt{1-}} & & \tilde{A}_4 \ar@{-}[dl]|*+{\txt{1-}} & &
S_4\!\!\times\! C_2 \ar@{-}[dr]|*+{\txt{2-}} & GL(2,3) \ar@{-}[d]|*+{\txt{2-}} & B \ar@{-}[dl]|*+{\txt{2-}} & \tilde{S}_4 \ar@{-}[dll]|*+{\txt{2-}} & &
A_5\!\!\times\! C_2 \ar@{-}[dr]|*+{\txt{1-}} & & \tilde{A}_5 \ar@{-}[dl]|*+{\txt{1-}} \\
& A_4 & & & & S_4 & & & & & A_5 }
\]
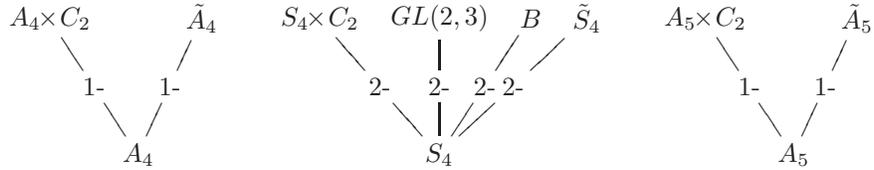
\centerline{Diagram \thediag : The double coverings of the rotary platonic groups.}
%

\section{Double coverings of some simple groups}\label{DCSG}

\subsection{Double coverings of projective linear groups over odd prime fields}

Let $p$ be an odd prime and $k=\frac{p+1}{2}$. Then
\[
P=\left\langle x,y\mid x^2,(xy)^3,(xy^4xy^k)^2y^p\right\rangle
\]
is a presentation of the projective special linear group $PSL(2,p)$ (see \cite{CR}). There are 2 presentation classes $[P_J]$
and for each choice of $j_2,j_3\in\{1,i\}$ we have:
\begin{itemize}
\item[1)] $P_{(1,j_2,j_3)}$ is a presentation of $PSL(2,p)\times C_2$.
\item[2)] $P_{(i,j_2,j_3)}$ is a presentation of $SL(2,p)$.
\end{itemize}

\[\addtocounter{diag}{1}
\xymatrix @R=36pt {
PSL(2,p)\!\!\times\! C_2 \ar@{-}[dr]|*+{\txt{1-}} & & SL(2,p) \ar@{-}[dl]|*+{\txt{1-}}\\
& PSL(2,p) }
\]
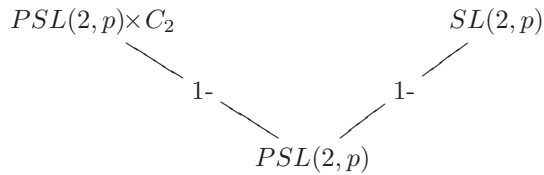
\centerline{Diagram \thediag : Double coverings of $PSL(2,p)$.}

\bigskip


\noindent Both double coverings are 1-strong, i.e. $Aut(SL(2,p))\cong Aut(PSL(2,p)\times
C_2)\cong Aut(PSL(2,p))$. Note that $A_4\cong PSL(2,3)$ and $A_5\cong PSL(2,5)$ were treated before as rotary platonic
groups.

\smallskip

The existence of only two double coverings is not a surprise. If fact, if $G$ is simple and non-abelian then $G$ is perfect ($G'=G$) and centerless.
Let $\hat G$ be a double covering of $G$. Then $\hat G$ is perfect or $\hat G'\lhd_2\hat G$. Since $G$ has no central
involution then the center $Z$ of $\hat G$ is $C_2$.

\smallskip

\noindent (i) The central extension $(\hat G,Z)$ is not irreducible. This means that $\hat G=NZ$ for some $N\lhd_2 \hat
G$. But then $N\cap Z=1$ and so $\hat G=N\times Z\cong N\times C_2\cong G\times C_2$. Reciprocally, if $\hat G=G\times
C_2$ then clearly $\hat G$ is not irreducible. Thus $(\hat G,Z)$ is not irreducible if and only if $\hat G=G\times
Z\cong G\times C_2$.

\smallskip

\noindent (ii) The extension $(\hat G,Z)$ is irreducible (that is, $\hat G$ is not a direct product $G\times C_2$).
Since $G$ is perfect, if $\hat G'\lhd_2\hat G$ then $\hat G'\cap Z=1$ and $\hat G=\hat G'\times Z\cong G\times Z$ is
not irreducible. Hence $\hat G$ must be also perfect. Consequently $\hat G'\cap Z=\hat G\cap Z=Z$ and as
$M(G)=Z\cap\hat G'$ (Theorem 9.9 \cite{Suzuki} (pg 251)) then $|M(G)|=|\hat G'\cap Z|=2$, where $M(G)$ is the
multiplier of $G$. Hence $(\hat G,Z)$ is a primitive central extension (by definition).
From Theorem 9.18(4) \cite{Suzuki} (pg 257) $\hat G$ is a representation group of $G$. Thus, if $(\hat G,Z)$ is an
irreducible double extension of $G$ then $\hat G$ is a representation group of $G$.

\medskip

For $G=PSL(2,q)$ with $q=p^f>4$ and $q\neq 9$,\, $SL(2,q)$ is the only representation group of $PSL(2,q)$ (Theorem
25.7, pg 646 of \cite{BH}). Hence $PSL(2,q)$ with $p\neq 2$ has only two double coverings $PSL(2,q)\times C_2$ and
$SL(2,q)$. If the field has characteristic 2 then $PSL(2,q)=SL(2,q)$ and hence $PSL(2,2^f)$ for $f>2$ has only one
double covering which is the direct product $PSL(2,2^f)\times C_2$.

\subsection{Other simple groups}

It is a known fact that simple groups are 2-generated \cite{AG}. Some simple groups $G$ have presentations
of simple type with only two odd relators. They have only one double covering, the direct product $G\times
C_2$. The unitary groups $U_3(3)$ and $U_3(4)$, the projective special linear group $PSL(3,3)$
and the Mathieu group $M_{11}$, as can be observed in Table 4, are such examples.
\begin{center}\addtocounter{tab}{1}
\begin{tabular}{llllr}
$G$ & Relators of $P$ & type of $P$ & $\delta$ & Ref.\\
\hline\\[-6pt]
$A_6$ & $a^4,\,b^5,\,abab^4abaBAB$ & simple & 1 & \cite{CHRR}\\
$A_7$ & $A^3(ab^2)^4,\, b^5,\,(b^2Aba)^2$ & simple & 2 & \cite{CRKMW}\\
$A_8$ & $(a^2b)^2,\,a^7,\,abAB^3AbAB^2$ & simple & 2 & \cite{CHRR}\\
$PSL(3,3)$ & $a^2B^3,\,BA(ba)^5(BA)^7(bA)^3(Ba)^2a^2B^2a(Ba)^2(ba)^2$ & simple & 1 & \cite{CR2}\\
$PSL(2,8)$ & $abAbaB,\, a^4(b^2ab)^2b$ & & 1 & \cite{CRKMW}\\
$PSL(2,16)$ & $a^3b^2A^2b^2,\, abAB^2AbaB$ & simple & 1 & \cite{CHRR}\\
$PSL(2,25)$ & $a^5,\, a^2b^2a^2b^2,\, ab^3a^4b^3aB$ & simple & 2 & \cite{CHRR}\\
$PSL(2,27)$ & $(ab)^2,\,a^7,\,a^2ba^6baBAb^4$ & simple & 2 & \cite{CHRR}\\
$PSL(2,32)$ & $abAbaB,\, (aba)^2b^3a^5b^3$ & & 1 & \cite{CRKMW}\\
$PSL(2,49)$ & $a^4,\, b^5,\, a^2b^4a^2BababaB$ & simple & 2 & \cite{CHRR}\\
$U_3(3)$ & $B^2ABa^3BA,\,b^2AB^2Ab^2aBa$ & simple & 1 & \cite{HNB}\\
$U_3(4)$ & $a^2B^3,\, a^3b(aB)^2(ab)^2aBab(aB)^4abaB(BA)^8b$ & simple & 1 & \cite{CR2}\\
$M_{11}$ & $aba^4b^3,\,babABaBabA$ & & 1 & \cite{K}\\
$M_{12}$ & $a^2B^3,\, (ab)^{10}[b,a]^6,\, \left( (ab)^4aBabaB\right)^3$ & & 2 & \cite{CR2}\\
$Sz(8)$ & $a^7,\,b^5,\, AB^2a^3ba^2B,\, BaBABa^5b^2a^2ba$ & simple & 4 & \cite{CHRR}\\
\hline
\end{tabular}\\[2pt]
Table \thetab: Presentations of some simple groups. Here $A=a^{-1}$ and $B=b^{-1}$.
\end{center}

\noindent As we saw before, the only double covering of $PSL(2,2^f)$, for $f\geq 3$, is
$PSL(2,2^f)\times C_2$. As it can be observed in Table 4, $PSL(2,8)$, $PSL(2,16)$ and $PSL(2,32)$
have deficiency zero. Do all $PSL(2,2^f)$ have a presentation of deficiency zero?

\medskip

Being $G$ simple, the double covering $G\times C_2$ is always 1-strong (Corollary \ref{SimpleStrong}). In the case of
$\delta=2$, besides $G\times C_2$ there is another double covering $\hat G\not\cong G\times C_2$ (Theorem
\ref{UniquePClass}) which is strong by Theorem \ref{strong}.

\vspace{50pt}

\begin{center}
\begin{tabular}{ccc}
Ana Breda &\quad Antonio Breda d'Azevedo &\quad Domenico Catalano\\
ambreda@mat.ua.pt &\quad breda@mat.ua.pt &\quad domenico@mat.ua.pt
\end{tabular}

\vspace{20pt}

Departamento de Matem\'atica\\
Universidade de Aveiro\\
Campus Universit\'ario de Santiago\\
PT-3810-193 Aveiro\\
Portugal
\end{center}

\end{document}